\documentclass[oneside,11pt]{article}
\usepackage{graphicx}

\begin{document}

\title{Genus two Smale-Williams solenoid attractors in
3-manifolds  \textit{} }

\author {Jiming Ma and Bin Yu  }

\maketitle

%\documentclass[oneside,12pt]{article}
%\textwidth=465truept \textheight=620truept \oddsidemargin=-1mm
%\renewcommand{\baselinestretch}{1.5}
%\usepackage{graphicx}
%\begin{document}

%\begin{center}

%\textbf{Genus two Smale-Williams solenoid attractors in
%3-manifolds}
%\end{center}

%\begin{center}
%\textbf{Jiming Ma and Bin Yu}

%\end{center}

\hspace*{-0.2cm} \textbf{Abstract:} Using alternating Heegaard
diagrams, we construct some 3-manifolds which admit diffeomorphisms
such that the non-wandering sets of the diffeomorphisms are composed
of Smale-Williams solenoid attractors and repellers, an interesting
example is the truncated-cube space. In addition, we prove that if
the nonwandering set of the diffeomorphism consists of genus two
Smale-Williams solenoids, then the Heegaard genus of the closed
manifold is at most two.

\hspace*{-0.2cm} \textbf{Keywords:}3-manifolds, Smale-Williams
solenoid attractors, alternating Heegaard diagram.

\hspace*{-0.2cm} \textbf{MR(2000)Subject Classification:} 57N10,
58K05, 37E99, 37D45.

 \vspace*{0.5cm}
\begin{bfseries}
1.Introduction
\end{bfseries}
\vspace*{0.5cm}

For a diffeomorphism  of a manifold $f:M \rightarrow M $, Smale
introduced the notion of hyperbolic structure on the non-wandering
set, $\Omega(f)$, of $f$. It is Smale's long range program to
classify  a Baire set of these diffeomorphisms, and $\Omega(f)$
plays a crucial role in this program. He also introduced solenoid
into dynamics in [S], in the literature this solenoid is called
Smale solenoid or pure solenoid.

To carry out this program, Williams defined 1-dimensional solenoid
in terms of 1-dimensional branched manifold, which is the
generalization of  Smale solenoid. There are two methods to define
Smale-Williams solenoid: the inverse limit of an expanding map on
branched manifold or the nested intersections of handlebodies.

Bothe studied the ambient structure of attractors in [B1], through
this work, we can see the two definitions above are equivalent.
Boju Jiang, Yi Ni and Shicheng Wang  studied the global question
in [JNW]. The question is, if a closed 3-manifold admit a
diffeomorphism $f$ such that the non-wandering set of $f$ consists
of two Smale solenoids, what we can say about the manifold $M$.
They proved that the manifold must be a lens space. Furthermore,
for any lens space, they can construct such a diffeomorphism. Our
previous paper [MY] also considered this question, we got all
Smale solenoids realized in a given lens space through an
inductive construction. Actually, part work of [JNW] and [MY] have
been studied in [B2].

A manifold $M$ admitting a diffeomorphism $f$ such that $\Omega(f)$
consists of two hyperbolic attractors presents a symmetry of the
manifold with certain stability. In the paper [JNW], the authors
noted that they believe many more 3-manifolds admit such symmetry if
we replace the Smale solenoid by its
generalization------Smale-Williams solenoid.

 In this paper we consider this question, more precisely, we have the following problem:
 is there any closed orientable three  manifold $M$ which  admits a
 diffeomorphism $f$ such that the non-wandering set of $f$,
 $\Omega(f)$,
 is composed of Smale-Williams
solenoid attractors and repellers?

In fact, for any 3-manifolds $M$, there is  a diffeomorphism $f$
such that $(M,f)$ has a Smale-Williams solenoid as one attractor,
so in the question, we must require that all of $\Omega(f)$
consists of Smale-Williams solenoids, and using standard arguments
in dynamics, in this case, there is exactly one attractor and one
repeller.

Gibbons studied this question in $S^{3}$ in [G]. He constructed many
such diffeomorphisms on $S^{3}$. Similar with the discussion in
[JNW], we focus on the question that which manifold admit such a
diffeomorphism. The main results of this paper are the following\\

\textbf{Theorem 3.4}:
 \begin{itshape} Let $M$ be a closed 3-manifold,
 and there is  $ f \in Diff(M)$ such that $\Omega(f)$ consists of
 genus two
  Smale-Williams
solenoids, then the Heegaard genus of $M$, $g(M)\leq 2$.
 \end{itshape}\\

\textbf{Theorem 4.5}:
\begin{itshape}
If a Heegaard splitting  $M=N_1\cup N_2$ of the closed orientable
3-manifolds $M$ is a genus two alternating Heegaard splitting,
then there is a diffeomorphism $f$, such that $\Omega(f)$ consists
of two Smale-Williams solenoids.
\end{itshape}\\

In fact, we give the first example that genus two Smale-Williams
solenoids can be realized globally in a Heegaard genus two closed
3-manifold. An interesting example is the rational homology sphere
whose fundamental group is the extended triangle group of order
48, i.e, the truncated-cube space, see [M].

Some notions in 3-dimensional manifolds theory and in dynamics will
be given in Section 2, for the definition of alternating Heegaard
splitting, see Section 4.

\vspace*{0.5cm}
\begin{bfseries}
2. Notions and facts  in 3-dimensional manifolds theory and in
dynamics
\end{bfseries}
\vspace*{0.5cm}

For fundamental facts about 3-manifolds see [H] and [J]. Let $G$ be
a finite graph in $R^3$, then a regular neighborhood $H$ of $G$ in
$R^3$ is called a\emph{ handlebody}, it is a 3-manifold with
boundary, the genus of its boundary is called the genus of $H$,
denoted by $g(H)$. Let $M$ be a closed orientable 3-manifold, if
there is a closed orientable surface $S$ in $M$ which separates $M$
into two handlebodies $H_1$ and $H_2$, then we say $M=H_1\cup_{S}
H_2$ is a\emph{ Heegaard splitting} of $M$, $S$ is called a
\emph{Heegaard surface}. Any closed orientable 3-manifold has
infinitely many  Heegaard splittings, and the minimum of the genus
of the Heegaard surfaces is called the Heegaard genus of the
3-manifold $M$.

A properly embedded 2-sided surface $F$ in a 3-manifold $M$ is
called an \emph{incompressible surface} if it is
$\pi_{1}$-injective, otherwise, it is a compressible surface.

The following theorems will be used in the paper:\\

\textbf{Haken Finiteness Theorem.}\emph{ Let $M$ be a compact orientable
3-manifold. Then the maximum number of pairwise disjoint,
non-parallel closed connected incompressible surfaces in $M$,
denoted by $h(M)$, is a finite integer $\geq 0$}.\\

\textbf{Papakyriakopoulos Loop Theorem.} \emph{Let $M$ be a
compact orientable 3-manifold and $S \subset M$ a closed
orientable surface. If the homomorphism $i_{*} :
\pi_1(S)\rightarrow \pi_{1}(M)$ induced by the embedding $i : S
\rightarrow M$ is not injective, then there is an embedded disk
$D\subset M$ such that $D\cap S = \partial D$
and $\partial D$ is an essential circle in $S$.}\\

We recall some facts about Smale-Williams
solenoid from the famous paper [W1].\\

\textbf{Definition 2.1}: A \emph{branched 1-manifold $L$}, is just
like a smooth 1-manifold, but there are two type of coordinate
neighborhoods are allowed. These are the real line $R$ and
$Y=\{(x,y)\in R^2:y=0$ or $y=\varphi(x)\}$. Here
$\varphi:R\longrightarrow R$ is a fixed $C^\infty$ function such
that $\varphi(x)=0$ for $x\leq 0$ and $\varphi(x)>0$ for $x> 0$.
The \emph{branch set $B$}, of $L$, is the set of all points of $L$
corresponding to $(0,0)\in Y$.  The first Betti number
$\beta_{1}(L)$ is called the \emph{genus }of the branched
1-manifold, and which is just the genus of the handlebody that $L$
induced(See Example 2.6).

Note that a branched 1-manifold $L$ has a tangent bundle $T(L)$,
and a differentiable map $f: L \rightarrow  L'$ between branched
1-manifolds induces a map $Df: T(L)\rightarrow T(L')$ of their
tangent bundles.

\begin{center}
\includegraphics[totalheight=4cm]{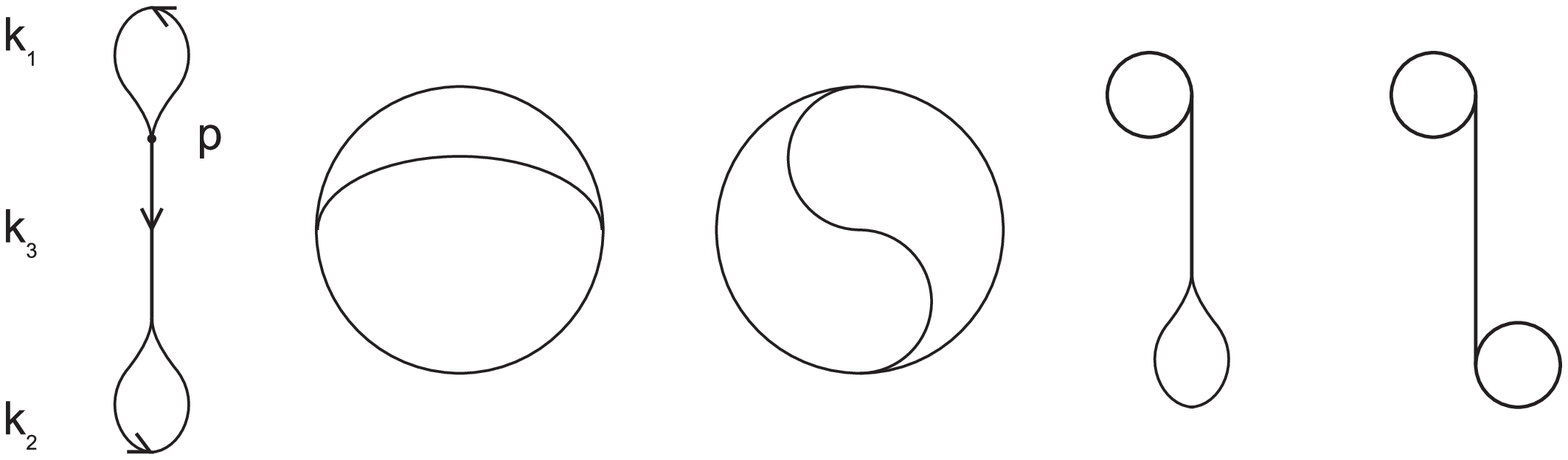}
\begin{center}
Figure 1
\end{center}
\end{center}

\textbf{Definition 2.2}: Let $L$ be a branched 1-manifold, a
$C^{r}$ immersion $g:L\longrightarrow L$ is called an
\emph{expansion map}, if there are $c>0,\lambda>1$, such that
$$\|(Dg)^{n}(v)\| \geq c \lambda^{n} \|v\|$$ $\forall n \in N$, $\forall v\in
T(L)$.\\

 \textbf{Definition 2.3}: Let $L$ be a branched 1-manifold with branch set $B$. We
 call
$g:L\longrightarrow L$ a
  \emph{Williams expansion map}, if:

 Axiom 1.  $g$ is an expansion map;

 Axiom 2.  $\Omega(g)=L$;

 Axiom 3.  Any point $p$ of $L$ has a neighborhood $U(p)$, such
  that $g(U(p))$ is an arc;

 Axiom 4.  There is a  finite set $A\subset L$, such that, $g(A\cup B)\subset
 A$.\\

 \textbf{Definition 2.4}: Let $\Sigma$ be the inverse limit of the sequence
 $$L  \stackrel{g}{\longleftarrow}  \ L\stackrel{g}{\longleftarrow}L\stackrel{g}{\longleftarrow}....$$
 where  $L$ is a branched 1-manifold, $g$ is a Williams expansion map on $L$.
 For a point a=($a_{0},a_{1},a_{2},...$) $\in \Sigma$, let
 $h^{-1}(a)=(a_{1},a_{2},a_{3},...)$, then $h:\Sigma \longrightarrow
 \Sigma$ is a homeomorphism. $\Sigma$ is called the \emph{Smale-Williams solenoid} with \emph{shift map} $h$, denoted it by
 $(\Sigma,h)$.\\

\textbf{Example 2.5}: See [W1], Figure 1 contains all of the
branched 1-manifolds with two branched points. Only the first two
allow immersions satisfing Axiom 2. Let $K$ be the first one and
define the Williams
 expansion map $g: K \rightarrow K$ on its oriented 1-cells by:
 $$K_1\longrightarrow K_3^{-1}K_1K_3$$
    $$  K_2\longrightarrow  K_3K_2K_3^{-1}$$
    $$  K_3\longrightarrow  K_2K_3^{-1}K_1$$

 It is easy to check that $g$ satisfies Axiom 1, 3, 4 of Definition 2.3.
 For Axiom 2, for some semi-conjugacy reason we only need to check that the induced
  symbolic dynamical system matrix $X$ is irreducible, that is, for
  all $1\leq i,j \leq dim(X)$, there exists $N(i,j)>0$ such that
  the $ij^{th}$ entry of $X^{N(i,j)}$ is positive, see [BH]. The induced matrix
  $X$ of $g$ is

 \begin{displaymath}
 \mathbf{X} =
 \left(\begin{array}{ccc}
  1&0&2\\
  0&1&2\\
  1&1&1
 \end{array}\right)
 \end{displaymath}
 obviously $X$ is irreducible.\\

\textbf{Example 2.6}: Figure 2(a) contains a genus two handlebody
$N$ with disk foliation and a self embedding $f$. In fact, we take
$N$ as a "neighborhood" of $K$, and there is a natural projection
$\pi$ from $N$ to $K$, and the embedding $f$ is induced by $g$ in
Example 2.5.
  We define $\Lambda=\bigcap_{n=1}^\infty  f^{n}(N)$, then $(\Lambda,f)$ is conjugate to $(\Sigma,h)$
via $T$ constructed below. We have the following commutative
diagram,
\[\begin{array}{lcl}
{N}& \stackrel{f}{\longrightarrow} & {N}\\
\downarrow \pi & & \downarrow \pi\\
{K}& \stackrel{g}{\longrightarrow} & {K}
\end{array}
\]
$\forall x \in \Lambda$, $T(x)$ is defined to be
$(\pi(x),\pi(f^{-1}(x)),...,\pi(f^{-n}(x)),...)$. As we shall see
in Proposition 2.7, $T$ is the conjugate map.

\begin{center}
\includegraphics[totalheight=6cm]{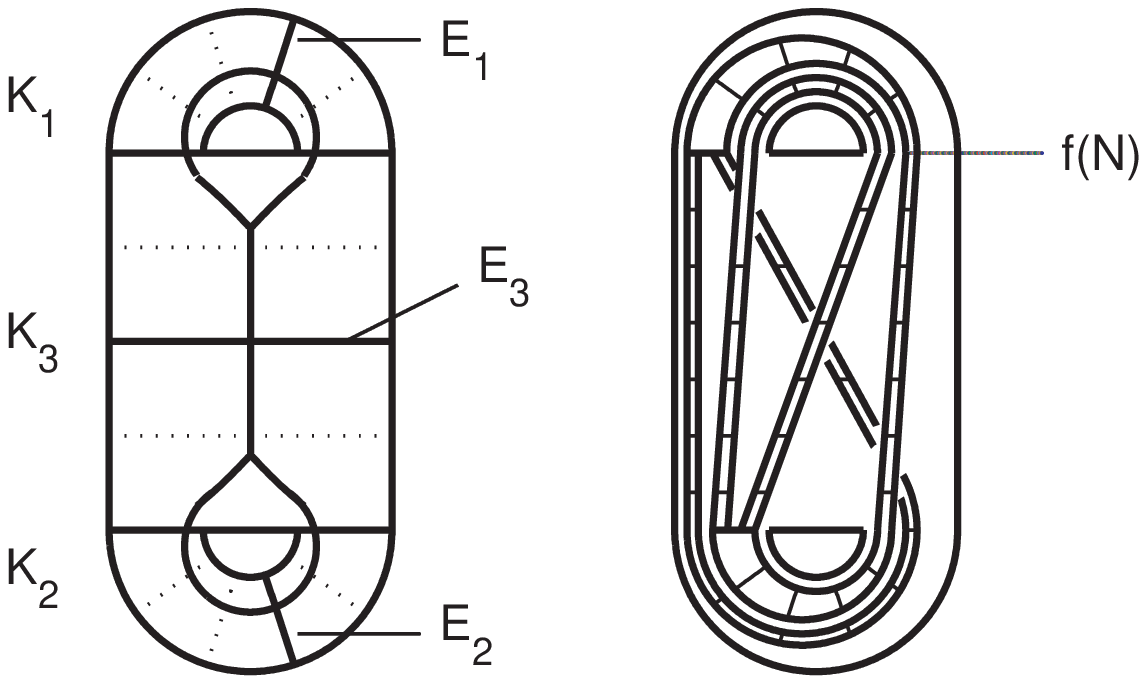}
\begin{center}
Figure 2(a)
\end{center}
\end{center}

\begin{center}
\includegraphics[totalheight=6cm]{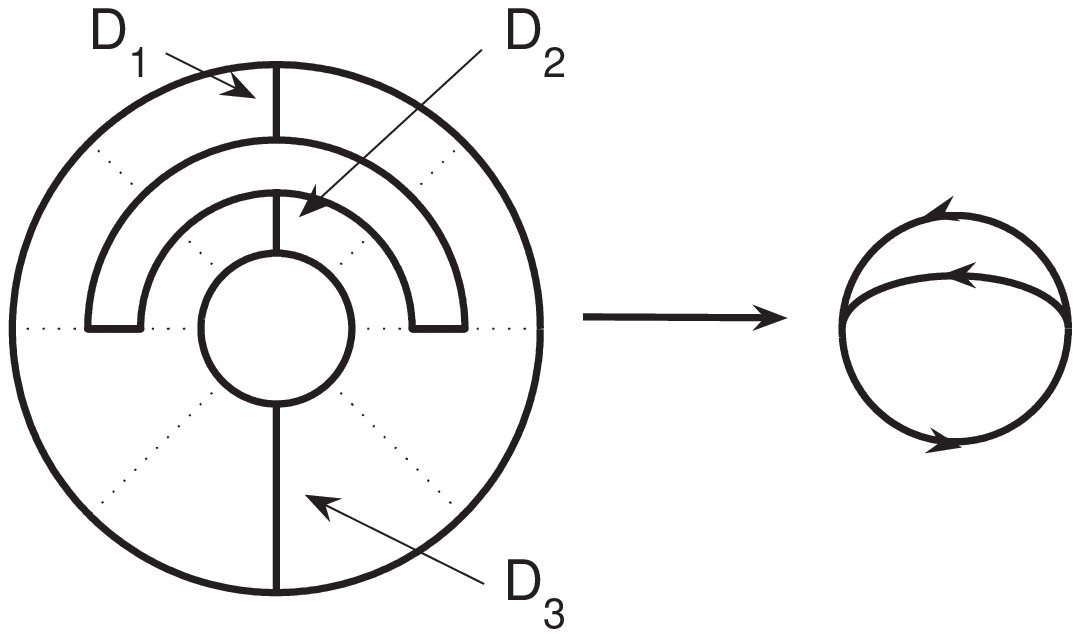}
\begin{center}
Figure 2(b)
\end{center}
\end{center}

Example 2.6 shows the way to construct Smale-Williams solenoid
locally in a geometrical way, the key is to deal with the branched
points,
we can get the following proposition:\\

\textbf{Proposition 2.7}:
\begin{itshape}
Let $L$ be a branched 1-manifold with Williams expansion map $g$
which induces Smale-Willams solenoid $\Sigma$ and shift map $h$.
Then the dynamical system $(\Sigma,h)$ has local three dimensional
model $(\Lambda,f)$, this means, there is a genus $\beta_{1}(L)$
handlebody $N$, an embedding $f:N\hookrightarrow N$, such that

  (1) There is a projection map $\pi:N\longrightarrow L$ which gives
a disc foliation structure of N by $\pi^{-1}(p)$, $\forall p\in L$.
$f$ preserves the foliation structure, furthermore, the area
satisfies
$$S_{f(\pi^{-1}(p))}/S_{\pi^{-1}(p)} \leq \epsilon$$ for some
$\epsilon>0$ small enough.

  (2) Define $\Lambda=\bigcap_{n=1}^\infty  f^{n}(N)$, we still
call $f|_{\Lambda}$ $f$, such that, $(\Sigma,h)$ is conjugate to
$(\Lambda,f)$.
\end{itshape}\\

\textbf{Proof}: This proof is based on the Example 2.6.

  \emph{1. Construct $(N,\pi,f)$.}

   Step a. Take a "neighborhood" of $L$, which is a
   handlebody $N$, and we foliated $N$ by disks as in Figure 2(a), so
   there is a natural
projection $\pi:N \rightarrow L$.

   Step b. Obviously there is a  local immersion map $i:
L\longrightarrow N$ such that $\pi\circ i=g$. Because of
$dim(N)=3=2dim(L)+1$, we can disturb $i(L)$ in $N$ such that we
get an embedding map $i_{1}: L\hookrightarrow N$ such that
$\pi\circ i_{1}=g$. Now we take a neighborhood of  $i_{1}(L)$, say
$N_{1}$,  such that $N_{1}$ is embedded into $N$. Just like step
a, it is easy to get an embedding map $f: N\hookrightarrow N$ such
that, $f(N)=N_{1}$ and satisfies condition (1) of Proposition 2.7.

\emph{2. Conjugation.}

   Let $\Sigma$ be the Smale-Williams solenoid defined by $(L,g)$, and  $h$ is the shift
map. Define $T: \Lambda \longrightarrow \Sigma$, $\forall x\in
\Lambda$, $T(x)= (\pi(x), \pi(f^{-1}(x)), ..., \pi(f^{-n}(x)),
...)$. Now it is easy to check that $T\circ f = h\circ T$ ,
$T\circ f(x)$=$(\pi\circ f(x)), \pi(x), \pi(f^{-1}(x)), ...)$,
$h\circ T(x)= (g\circ \pi(x),\pi(x), \pi(f^{-1}(x)), ... )$. Since
$\pi\circ f = g\circ \pi$, we get $T\circ f = h\circ T$. On the
other hand, T has inverse map $T^{-1}: (x_{0}, x_{1}, x_{2}, ...
)\longrightarrow \pi^{-1}(x_{0})\cap f(\pi^{-1}(x_{1}))\cap
f(\pi^{-2}(x_{2})) \cap ... $, this map is well-defined by
$S_{f(\pi^{-1}(p))}/S_{\pi^{-1}(p)} \leq \epsilon$, it is easy to
check it is the inverse map of T. So $(\Sigma,h)$ is conjugate to
$(\Lambda,f)$ by $T^{-1}$. Obviously, $(N,f)$ constructed above is
an attractor
model. Q.E.D.\\

This type of attractors have important meaning in the study of
attractors, see Williams [W2].

 \vspace*{0.5cm}
\begin{bfseries}
3. Restriction on the Heegaard genus
\end{bfseries}
\vspace*{0.5cm}

We term the first two branched 1-manifolds \emph{type I} and
\emph{type II }respectively, and discuss their Williams expansion
map using their geometrical model introduced above. We call this
type of
attractor \emph{type I(type II) Smale-Williams solenoid attractor}.\\

\textbf{Definition 3.1}: Let $N$ be a genus two handlebody in a
3-manifold $M$, if there is  $ f \in Diff(M)$
 such that $f|_{N}$ is conjugate to a local model of type I(type II) Smale-Williams solenoid attractor,
then we call $(M,f)$ has \emph{type I(type II) Smale-Williams
solenoid
 attractor}.\\

\textbf{Lemma 3.2}: \emph{Let $(N, f)$ be a local model of a genus
two Smale-Williams solenoid attractor,
 then there is no properly embedded disk $D$ in $N$ such that  $f(N)\cap
 D=\emptyset$}.\\

\textbf{Proof:} We prove the lemma only for type I Smale-Williams
solenoid attractor, and the proof for the case of type II is
similar.

If there is an essential disk $D$ such that $f(N)\cap D=\emptyset$,
then there is a solid torus $V\subset N$ such that $Im (f)\subset
V$, this induces $Im(f_{\ast}) \subset \pi_{1}(V)\cong \
\textbf{Z}$, so $Im(f_{\ast})$ is abelian.

Let $K$ be a type I branched 1-manifold as Figure 1 shows and  $g$
be the induced map of $f$ on $K$, which is a train track map, see
[BH].

% We claim that $f_{\ast}$ is injective, here $f_{\ast}\in
%Hom(\pi_{1}(N))$ is induced by $f$. This is equivalent to $
%g_{\ast}$ is injective, $g_{\ast}\in Hom(\pi_{1}(K))$ is induced by
%$g$.

We choose a base point $P$ for $\pi_{1}(K)$, $x=[K_1]$, $y=[K_3 K_2
K^{-1}_{3}]$, $\pi_{1}(K)=<x>*<y>$ is a free group of rank two(see
Figure 1). Since $Im(g_{\ast})$ is  an abelian group, we have that
$g_{\ast}(xyx^{-1}y^{-1})=1$. $xyx^{-1}y^{-1}=[K_1 K_3 K_2
K^{-1}_{3} K^{-1}_1 K_3 K^{-1}_2 K^{-1}_{3}]$. Since $K_1 K_3 K_2
K^{-1}_{3} K^{-1}_1 K_3 K^{-1}_2 K^{-1}_{3}$ is a legal path in $K$
 and $g$ is a train track map, we have $g(K_1 K_3 K_2
K^{-1}_{3} K^{-1}_1 K_3 K^{-1}_2 K^{-1}_{3})$ is also a legal path,
so it can not be a homotopic trivial path  in $K$(See [BH]). This
means $g_{\ast}(xyx^{-1}y^{-1})\neq 1$, it is a contradiction.
Q.E.D.
\\

\textbf{Proposition 3.3}: \emph{Let $M$ be a closed 3-manifold,
and $N$ is a genus two handlebody in $M$. If there is $ f\in
Diff(M)$
 such that $(N,f|_{N})$ is conjugate to a local genus two Smale-Williams solenoid, then $\partial N$ is
 compressible in $\overline{M-N}$}.\\

 \textbf{Proof}: Suppose $\partial N$ is incompressible in $\overline{M-N}$. Let $m$ be
  the  Haken number of $\overline{M-N}$, denoted by
 $h(\overline{M-N})= m$,
  $S_{1}$, $S_{2}$, ..., $ S_{m-1}$, $\partial N $
 are mutually disjoint nonparallel  incompressible surfaces in
 $\overline{M-N}$. Since $f$ is a diffeomorphism  from
 $\overline{M-N}$ to $\overline{M-f(N)}$, so $h(\overline{M-f(N)}) = m$, and $\partial f(N)$ is incompressible in
 $\overline{M-f(N)}$. If $\partial N$ is compressible in
 $\overline{M-f(N)}$, then $\partial N$ is compressible in
 $\overline{N-f(N)}$ , which contradicts to Lemma 3.2.
Then, by standard arguments in 3-manifold topology, $S_{1}, S_{2},
..., S_{m-1}$ are incompressible surfaces in
 $\overline{M-f(N)}$.
And then,  $S_{1},
 S_{2}, ..., S_{m-1}$, $\partial N$, $\partial f(N)$ are mutually disjoint nonparallel
 incompressible surfaces in $\overline{M-f(N)}$, so
 $h (\overline{M-f(N)})\geq m+1$, which contradicts to
 $h(\overline{M-f(N)}) = m$.
   Q.E.D.\\

 \textbf{Theorem 3.4}:
 \begin{itshape}Let $M$ be a closed 3-manifold,
 and there is  $ f \in Diff(M)$ such that $\Omega(f)$ consists of
 genus two Smale-Williams solenoids, then the Heegaard genus of $M$, $g(M)\leq 2$.
 \end{itshape}\\

\textbf{Proof}:
 If the nonwandering set $\Omega(f)$ consists of genus two
 Smale-Williams solenoids, then standard arguments in dynamics
 theory shows
  $\Omega(f)=\Lambda_{1}\sqcup \Lambda_{2}$
where $\Lambda_{1}$ is an attractor, $\Lambda_{2}$ is a repeller
(See [JNW]). In addition, $\Lambda_{1}$($\Lambda_{2}$) is realized
by genus two handlebody $N_{1}(N_{2})\subset M$,
$\Lambda_{1}=\cap_{n\geq 0}f^{n}(N_{1})$ and
$\Lambda_{2}=\cap_{n\leq 0}f^{n}(N_{2})$. There are $m,k\in
Z^{+}$, such that $\partial f^{m}(N_{2})\subset
\overline{N_{1}-f^{k}(N_{1})}\subset N_{1}$, so if necessary, we
let
 $f^{m}(N_{2})$ be the new $N_{2}$, then
we have $\partial N_{2}\subset \overline{N_{1}-f^{k}(N_{1})}\subset
N_{1}$. By Proposition 3.3, $\partial N_{2}$ is compressible in
$\overline{M-N_{2}}=\overline{N_{1}-N_{2}}$. Let $c$ be an essential
simple closed curve in $\partial N_{2}$ which bounds a disk $D$ in
$\overline{N_{1}-N_{2}}$. Adding a neighborhood of $D$ in
$N_{1}-N_{2}$ to $N_{2}$ (2-handle addition along $c$), and denote
the resulting manifold by $ N_{2}^{*}$:

 \emph{Case 1.} If $c$ is non-separating in $\partial N_{2}$,  $ \partial N_{2}^{*}=S$ is a torus. In this case $N_{1}$ is divided into two parts $W_{1}$ and $W_{2}$ by
 $S$, $\partial W_{2}=S\sqcup \partial N_{1}$ and  $\partial W_{1}=S$:

 \emph{Subcase 1.1.} $S$ is compressible in $W_{1}$, so $W_{1}$ is a solid torus, and $\overline{N_{1}-N_{2}}$ is
 a genus two handlebody.  Hence $g(M)\leq 2$.

 \emph{Subcase 1.2.} $S$ is incompressible in $W_{1}$, since $N_{1}$ is a handlebody, $S$
 must be compressible in $W_{2}$. Compressing $S$ in  $W_{2}$,
 we get a $2$-sphere $P$ in $W_2$, so $P$ is also in $ N_{2}^{*}$.
  $P$ separates $W_2$ into $A$ and $B$, where $S \subset
 \partial A$ and $\partial N_1 \subset \partial B$. And we have $N_{2}^{*}= A^{*} \sharp
 B^{*}$, where  $A^{*}$ is obtained from $A$ by capping off $P$, and
 $B^{*}$ is obtained from $B\cup \overline{N_2-N_1}$ by capping
 off $P$. $W_1\cup A$ is in the handlebody $N_1$, so $W_1\cup A$ is a $3$-ball, $\pi_1(A^{*})$ is nontrivial. And
  the rank of $\pi_1( N_{2}^{*})$ is at most two, we get that the rank of
  $\pi_1 (B^{*})$ is at most one. So $ B^{*}$ has a genus one
  Heegaard splitting by the fact that the Heegaard genus of  $
  N_{2}^{*}$ is at most two. So  $M=B^{*}$ has Heegaard genus at
  most one.

 \emph{Case 2.}  If $c$ is separating in $\partial N_{2}$, $ \partial N_{2}^{*}=S$ is composed of two tori $S_{1}$ and $S_{2}$, so
 there are two manifolds $W_1$ and $W_2$ in $\overline{N_1-N_2}$
 with $\partial W_{i}= S_{i}$:

 \emph{Subcase 2.1.} One of $W_{1}$ and $W_{2}$ is a solid torus.
 Then there is a nonseparating simple closed curve in $\partial
 N_2$ which bound a disk in $\overline{N_1-N_2}$, and turn to Case 1.

\emph{Subcase 2.2.} $S_{i}$ is incompressible in $W_{i}$, so $S_{i}$
must be compressible in $\overline{N_1- W_{i}}$. Let $E$ be a
compressible disk of $S_{1}$ in $\overline{N_1- W_{1}}$, we can
assume $E\cap D=\emptyset$, $E\cap S_2$ is a set of simple closed
curves,
 and $|E\cap S_2|$ is minimal along all such compressible disk of $S_1$. If $E\cap S_2\neq \emptyset$
 , take an innermost disk, say $\Delta$, of $E$, which is a compressible disk of $S_2$ in
 $\overline{N_1- W_2}$, which is also disjoint from $D$. Compressing $S_2$ along $\Delta$, we get a $2$-sphere $P_2$,
  which bounds a $3$-ball
 $B_2$ in $N_1$, and $W_2 \subset B_2$. From this $3$-ball, and the
  compressibility of $S_1$ in $\overline {N_1-W_1}$, we can get a
 $3$-ball $B_1$, which is disjoint from $D$ and $B_2$, and contains $W_1$, $\partial B_1=P_1$. Connecting
  $P_1$ and  $P_2$ by $D$, we get a separating $2$-sphere $P$ in $ \overline{N_1 \cap N_2}$, which
  separates $M$ into two components,
  each component is a $B^3$ since it is contained in a handlebody. So $M=S^3$.
  Q.E.D\\

 \vspace*{0.5cm}
\begin{bfseries}
4. Globally geometric realization of Smale-Williams solenoid type
attractors in 3-manifolds
\end{bfseries}
\vspace*{0.5cm}

\textbf{Definition 4.1}: Let $l_1$, $l_2$ be two subarcs of a
branched 1-manifold $L$, $e$ an arc such that $e \cap (l_1\cup
l_2) =\partial e$, $E=I\times I=[-1,1]\times [-1,1]$ is a band
with the core $\{0\}\times [-1,1]=e$, $E \cap (l_1\cup l_2)
=[-1,1]\times \{-1,1\}$. Then along $e$, we get two new subarcs
$l_3$, $l_4$ of a new branched 1-manifold $L'$, $l_3$ and $l_4$
has just one crossing, this process is said to be a \emph{band
move} along $E$. And note that there are two band moves along $E$,
see Figure 3.

\begin{center}
\includegraphics[totalheight=3cm]{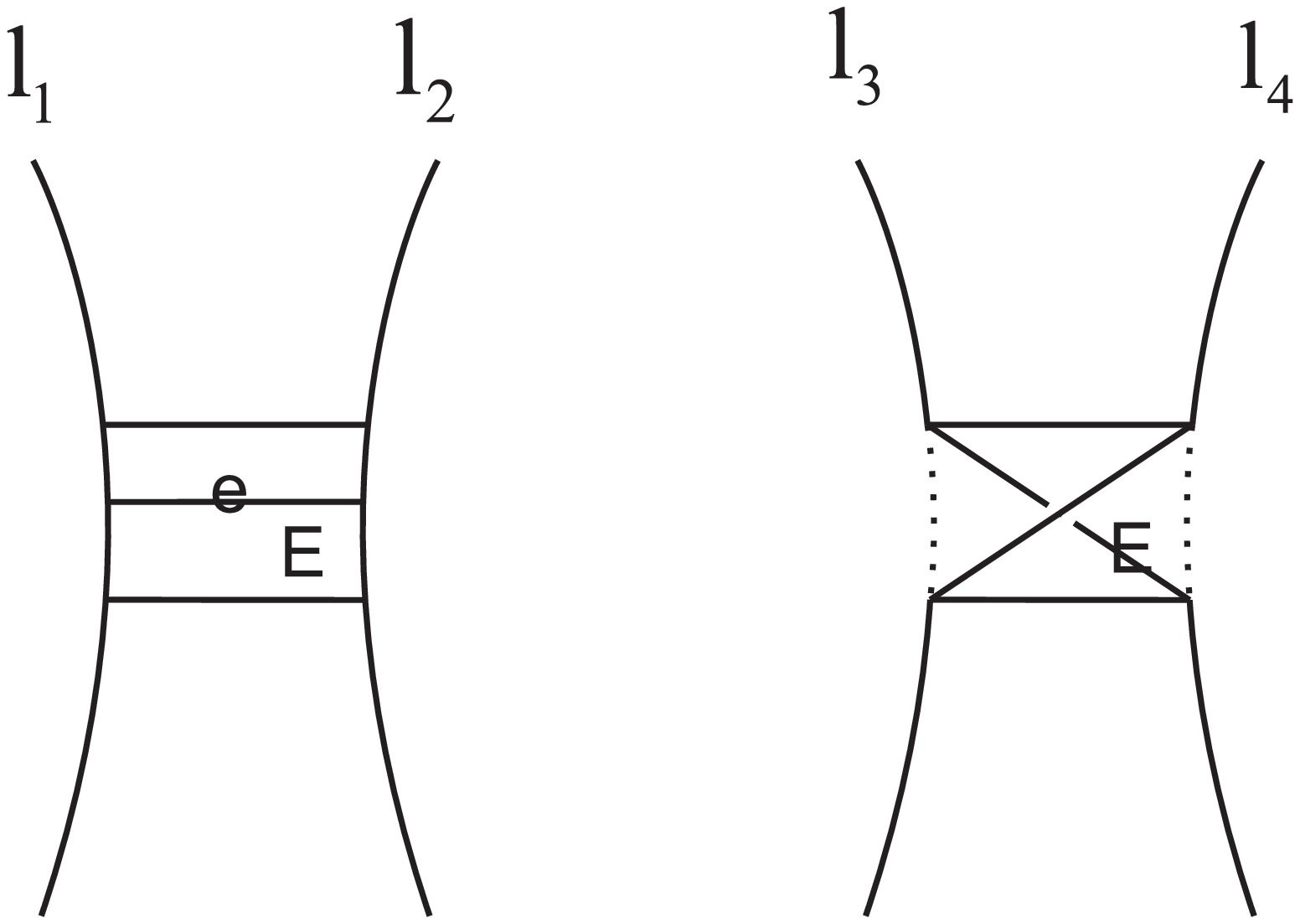}
\begin{center}
Figure 3
\end{center}
\end{center}

\textbf{Definition 4.2(Alternating Heegaard splitting of type I)}:
As in Figure 2(a), let $N$ be a handlebody, $K$ be a branched
1-manifold of type I, $\pi: N\rightarrow K$ be the natural
projection,
 then there are three disks in $N$, say $E_1$, $E_2$
 and $E_3$, which separate $N$ into two 3-balls, and $E_3$ is separating in $N$,  the image of $E_i$ in $K_i$ is an interior point in
 $K_i$. Let $c$ be a
  simple closed curve in $\partial N$, we say
$c$ is
 alternating with respect to $(E_1,E_2,E_3)$, if the intersection points of $E_{i}$ occur in $c$
alternatively about $E_3$, that is, along $c$,
  we see  $E_3$, $E_{i_{1}}$, $E_3$, $E_{i_{2}}$, $E_3$, $E_{i_{3}}$
  $E_3$...., where $i_{j}=1$ or $2$ and $i_{j}\cdot i_{j+1}=2$ to all $j$.

   Let
   $M=N_1\cup N_2$ be a genus two Heegaard splitting, if there
  are three disks $D_1$, $D_2$, $D_3$ separate $N_1$ into two
  3-ball, where $D_3$ is separating in $N_1$, and there
  are three disks $C_1$, $C_2$, $C_3$ separate $N_2$ into two
  3-ball, where $C_3$ is separating in $N_2$, moreover $\partial
  C_i$ is alternating with respect to  $(D_1,D_2,D_3)$, and $\partial
  D_i$ is alternating with respect to  $(C_1,C_2,C_3)$, then we
  say the Heegaard splitting is  \emph{alternating of type I}.

\textbf{Definition 4.3(Alternating Heegaard splitting of type II)}:
As in Figure 2(b), let $N$ be a handlebody, $L$ be a branched
1-manifold of type II,  $\pi: N\rightarrow L$ be the natural
projection, then there are three disks in $N$, say $E_1$, $E_2$
 and $E_3$, which separate $N$ into two 3-balls, none of $E_i$ is separating in $N$, the image of $E_i$ in $L_i$ is an interior point in
 $L_i$. Let $c$ be a
  simple closed curve in $\partial N$, we say
$c$ is
 alternating with respect to $(E_1,E_2,E_3)$, if the intersection points occur in $c$ alternatively, that is, along $c$,
 we see  $E_3$, $E_{i_{1}}$, $E_3$,
$E_{i_{2}}$, $E_3$, $E_{i_{3}}$, $E_3$...., where $i_{j}=1$ or
$2$.

   Let
   $M=N_1\cup N_2$ be a genus two Heegaard splitting, if there
  are three disks $D_1$, $D_2$, $D_3$ separate $N_1$ into two
  3-balls, none of $D_i$ is separating in $N_{1}$, and there
  are three disks $C_1$, $C_2$, $C_3$ separate $N_2$ into two
  3-balls, none of $C_i$ is separating in $N_{2}$. Moreover $\partial
  C_i$ is alternating with respect to  $(D_1,D_2,D_3)$, and $\partial
  D_i$ is alternating with respect to  $(C_1,C_2,C_3)$, then we
  say the Heegaard splitting is  \emph{alternating of type II}.\\

\textbf{Proposition 4.3}:
\begin{itshape}
If there is a  genus two handlebody $N$ in $M$, and an alternating
simple closed curve $c$ in $\partial N$ which bounds a disk in
$\overline{M-N}$,
 then there is $ f\in Diff(M)$, such that $\Omega(f)$ contains a genus two Smale-Williams
solenoid attractor.
\end{itshape}\\

\textbf{Proof}:  The proof is an explicit construction, we construct
the diffeomorphism in the type I case, the type II case is similar.

 We choose three parallel curves $c_{1}$, $c_{2}$, $c_{3}$ in $\partial
 N$ which are parallel to $c$. The branched 1-manifold $J$ is a spine of $N$, which
 is composed of oriented 1-cells: $J_{1}$, $J_{2}$, $J_{3}$ as
Figure
 4
shows. Note that $J$ induces a disk foliation of $N$, so
$\pi:N\rightarrow J$ is the projection map. We do the following
operations to $J$:

 \emph{Operation 1.} As indicated in Figure(4-1)$\rightarrow$
 Figure(4-2). We take a subarc $J_{1,1}$ of $J_1$, half-twist and move it
 toward $\partial N$ and identify it with a subarc $c_{1,1}$ of
 $c_{1}$, just like Figure (4-2) shows, this process is a band move in Definition 4.1.
 Since $c_{1}$ bounds a
 disk in $\overline{M-N}$, we can push $c_{1,1}$ across the disk
 to the arc of the subarc $c_{1,2}$, where $c_{1,2}=\overline{c_{1}-c_{1,1}}$.
  We do the same surgeries  to $J_{2}$, $J_{3}$. In the end, we get a
  branched 1-manifold
 $J^*$ which is isotopic to $J$ in $M$.

\begin{center}
\includegraphics[totalheight=10cm]{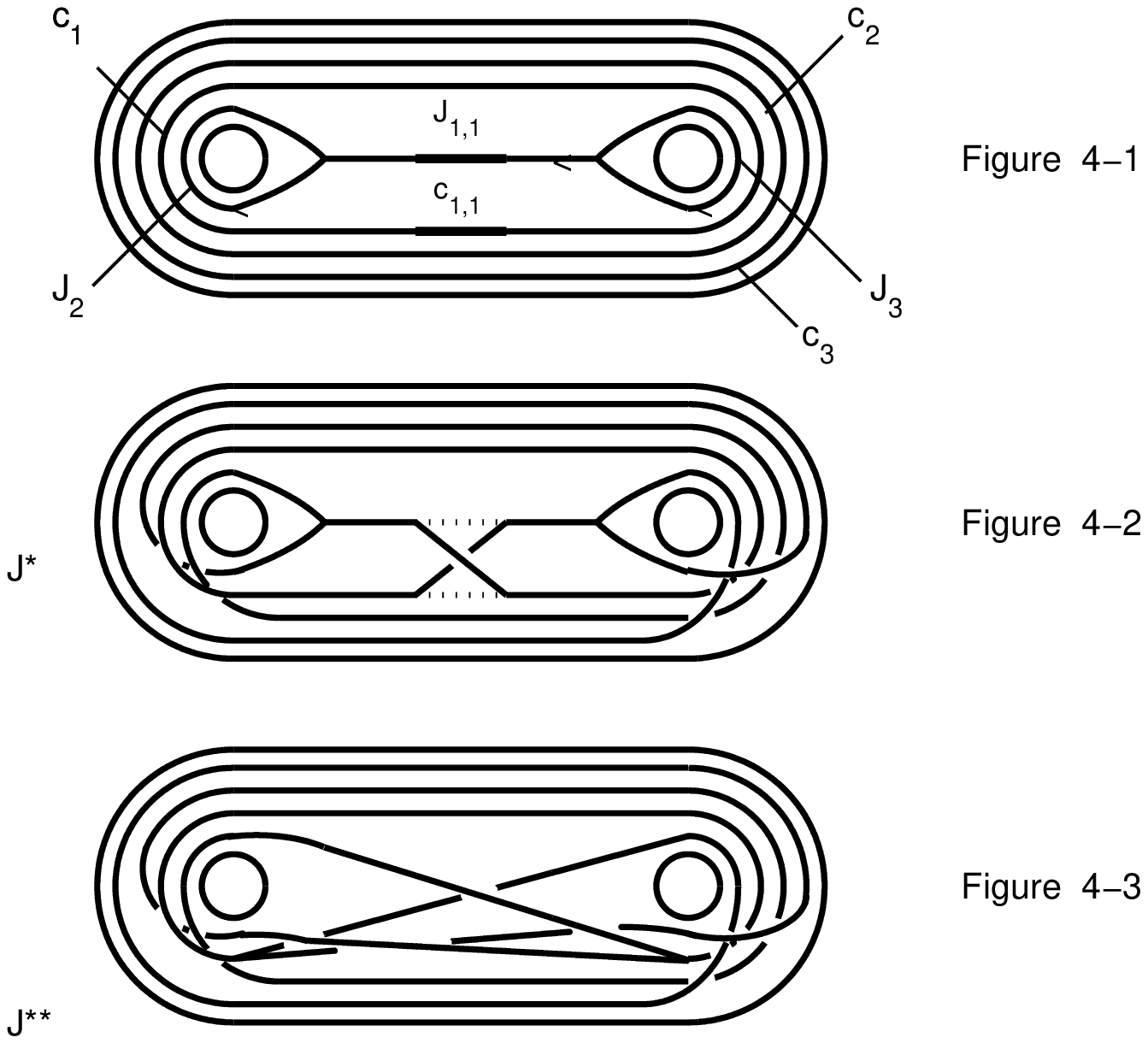}
\begin{center}
Figure 4
\end{center}
\end{center}

\begin{center}
\includegraphics[totalheight=1.5cm]{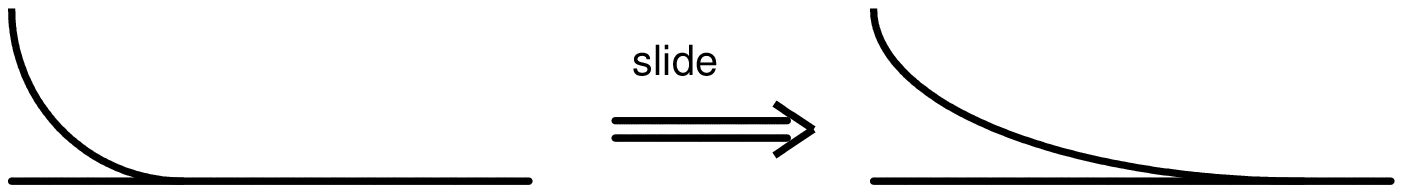}
\begin{center}
Figure 5
\end{center}
\end{center}

 \emph{Operation 2.} Note that the composite map $J \longrightarrow J^* \longrightarrow
 J$ is not a  Williams expansion map, where $J^* \longrightarrow
 J$ is induced by $\pi$.
So as indicated in Figure(4-2)$\rightarrow$ Figure(4-3). Slide one
end of $J_{2}^*$ along a subarc of $J_{1}^*$, also one end of
$J_{3}^*$ along another subarc of $J_{1}^*$, locally as Figure 5
shows. We get a new branched 1-manifold $J^{**}$
 which is isotopic to $J$ in $M$ obviously.

 So we can choose an $f\in Diff(M)$ which is isotopic to the identity such
 that:

 (1) $f:N\hookrightarrow N$;

 (2) $f(J)=J^{**}$;

 (3) $f$ preserve the disk foliation structure of
 $N$. For every leaf $D$ of the foliation, the area satisfies  $S_{f(D)}/S_{D} \leq \epsilon$ for some
$\epsilon>0$ small enough.

 (4) Let $g=\pi \circ f$, we get the following diagram,
 \[\begin{array}{lcl}
{N}& \stackrel{f}{\longrightarrow} & {N}\\
\downarrow \pi & & \downarrow \pi\\
{J}& \stackrel{g}{\longrightarrow} & {J}
\end{array}
\]
 $g$ is linear on every edge of $J$.

\emph{Claim:} $g$ is a Williams expansion map of $J$.

 \textit{Proof of the Claim.} We check Axiom 1, ... , Axiom 4 of Definition 2.3 one by
one.

 \emph{Axiom 1.} From the construction of $f:N\rightarrow
N$, for example in Figure 4, we know, $g$: $J_1 \rightarrow J_2
J^{-1}_1 J_3$,
   $J_3 \rightarrow J^{-1}_1 J_3 J_1 J_2 J^{-1}_1 J_3$, and
   $J_2 \rightarrow J_1 J_2 J^{-1}_1 J_3 J_1 J_2  J^{-1}_1$.
  Since $f$ expands wholly and $c$ has the alternating property then
  $g$ is smooth immersion for every local smooth arc of $J$. Since
  there are two branched points, $g$ is a Williams expansion map.

\emph{Axiom 2.} Since the matrix of the symbol dynamical system
induced by $g:J \rightarrow J$ is irreducible, actually nowhere is
zero in the matrix through the check of Axiom 1. Thus any point of
$J$ is a nonwandering point of $(J,g)$.

\emph{ Axiom 3 and Axiom 4.} These are obviously. So the Claim
follows.

By the Claim, we know $f\in Diff(M)$ and $\Omega(f)$
 contains a genus two Smale-Williams solenoid attractor. Q.E.D.\\

\textbf{Example 4.4($RP^{3}$)}:  Figure 6 is an alternating genus
2 Heegaard diagram of $M=RP^{3}=N_1\cup N_2$. The left figure
depicts the diagram seen from out of $N_1$, $d_i$ bounds a disk
$D_{i}$ in $N_1$, $c_i$ bounds a disk $C_{i}$ in $N_2$.
 The right figure depicts the
diagram seen from out of $N_2$,  it comes from the left diagram
 by Dehn twists $D_{e_3} D_{e_2} D_{e_2} D_{e_1}$ in Figure 7 and then follows by a mirror symmetry.

\begin{center}
\includegraphics[totalheight=8cm]{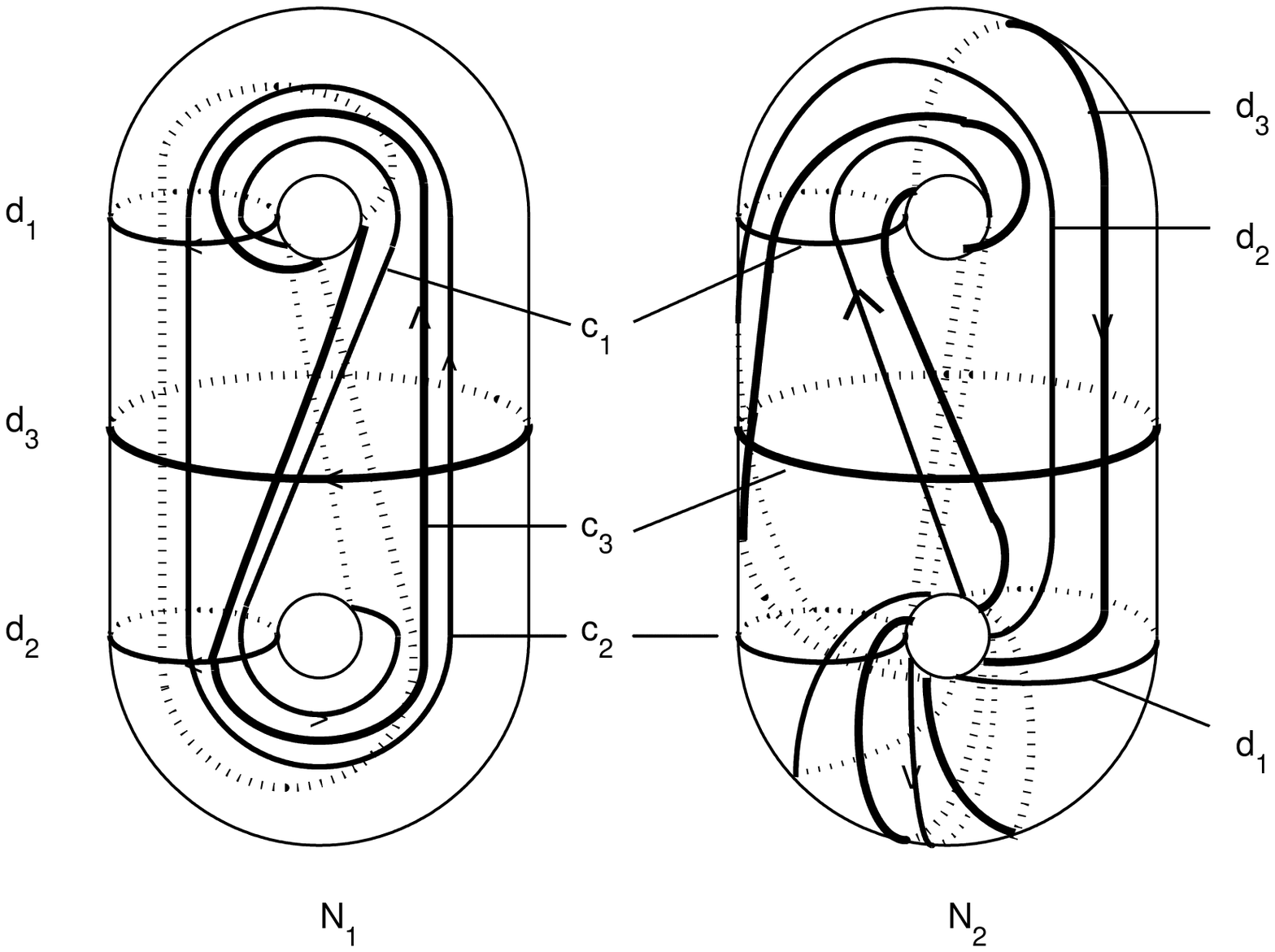}
\begin{center}
Figure 6
\end{center}
\end{center}

 And now,
 by the diagram, we  construct a self-diffeomorphism $f$ of $RP^{3}$
 such that
 $\Omega(f)$ consists  of two type I Smale-Williams solenoids.

\begin{center}
\includegraphics[totalheight=5cm]{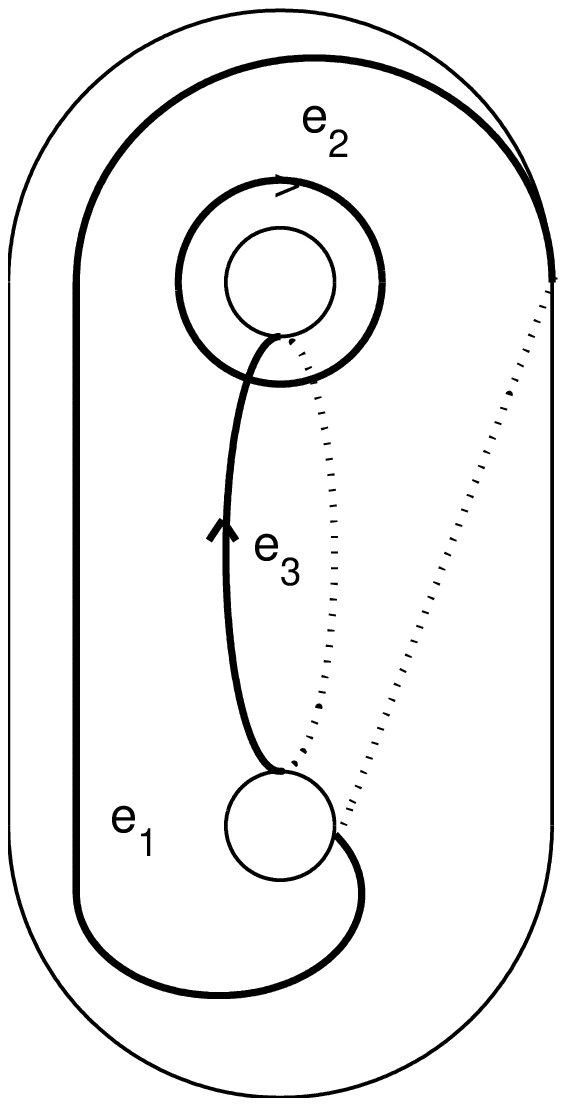}
\begin{center}
Figure 7
\end{center}
\end{center}

  Let $J$ be the natural spine of $N_1$ with respect to the three disks
  $D_1$, $D_2$ and $D_3$ bounded by $d_1$, $d_2$ and $d_3$,
  and $L$ be
the spine of $N_2$ with respect to the three disks $C_1$, $C_2$
and $C_3$ bounded by $c_1$, $c_2$ and $c_3$.

  Let $J\cap D_{i}=P_{i}$, $L\cap C_{i}=O_{i}$ and we fix a point $Q_{i}\in c_{i}\cap
  d_{i}$. We connect $P_{i}$ with $Q_{i}$ by an arc $v_i$ in $D_i$, and
  connect $Q_{i}$ with $O_{i}$ by an arc $w_i$ in $C_i$.

   We isotopy $J$ to $J^{*}$ in $M$ by performing three band moves along $v_{1}$,
$v_{2}$ and $v_{3}$, similarly we can isotopy $L$ to $L^{*}$ by
performing three band moves along  $w_{1}$, $w_{2}$ and $w_{3}$. By
Figure 8, we can see that $L^{*}\sqcup J$ is isotopic to
$J^{*}\sqcup L$: after three local half twist surgeries and pushing
moves to $J^{*}\sqcup L$ (see Figure 8-2), we get $L^{*}\sqcup
J$(see Figure 8-3), all the surgery can be regarded as appearing in
three mutually disjoint 3-balls $N(v_{i}\sqcup w_{i})$ in $M$(see
Figure 9), so this progress is an isotopy move, hence $L^{*}\sqcup
J$ is
isotopic to $J^{*}\sqcup L$. This process is similar with [JNW].\\

\begin{center}
\includegraphics[totalheight=2.5cm]{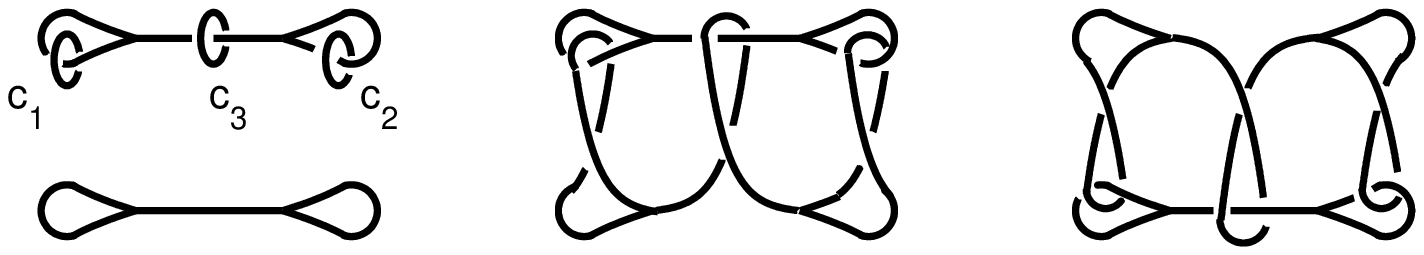}

\begin{center}
Figure 8-1 ~~~~~~~~~~~~~~~~~~~~Figure 8-2
~~~~~~~~~~~~~~~~~~~Figure 8-3
\end{center}

\end{center}

\begin{center}
\includegraphics[totalheight=5cm]{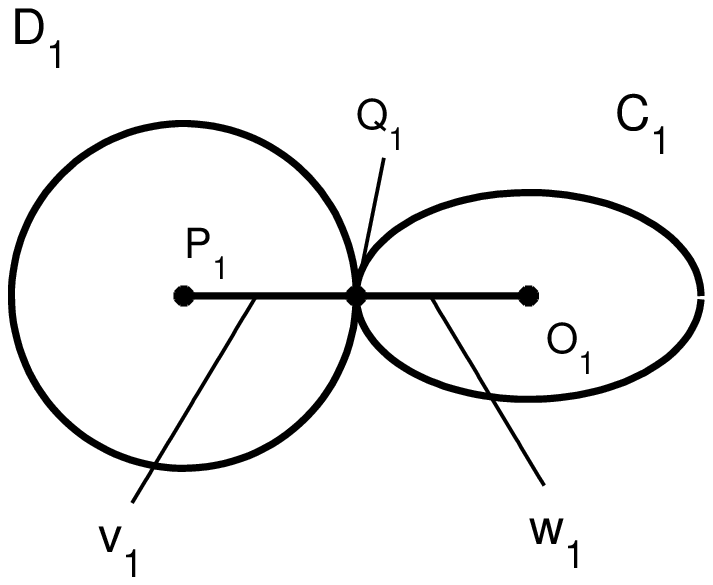}
\begin{center}
Figure 9
\end{center}
\end{center}

\begin{center}
\includegraphics[totalheight=6cm]{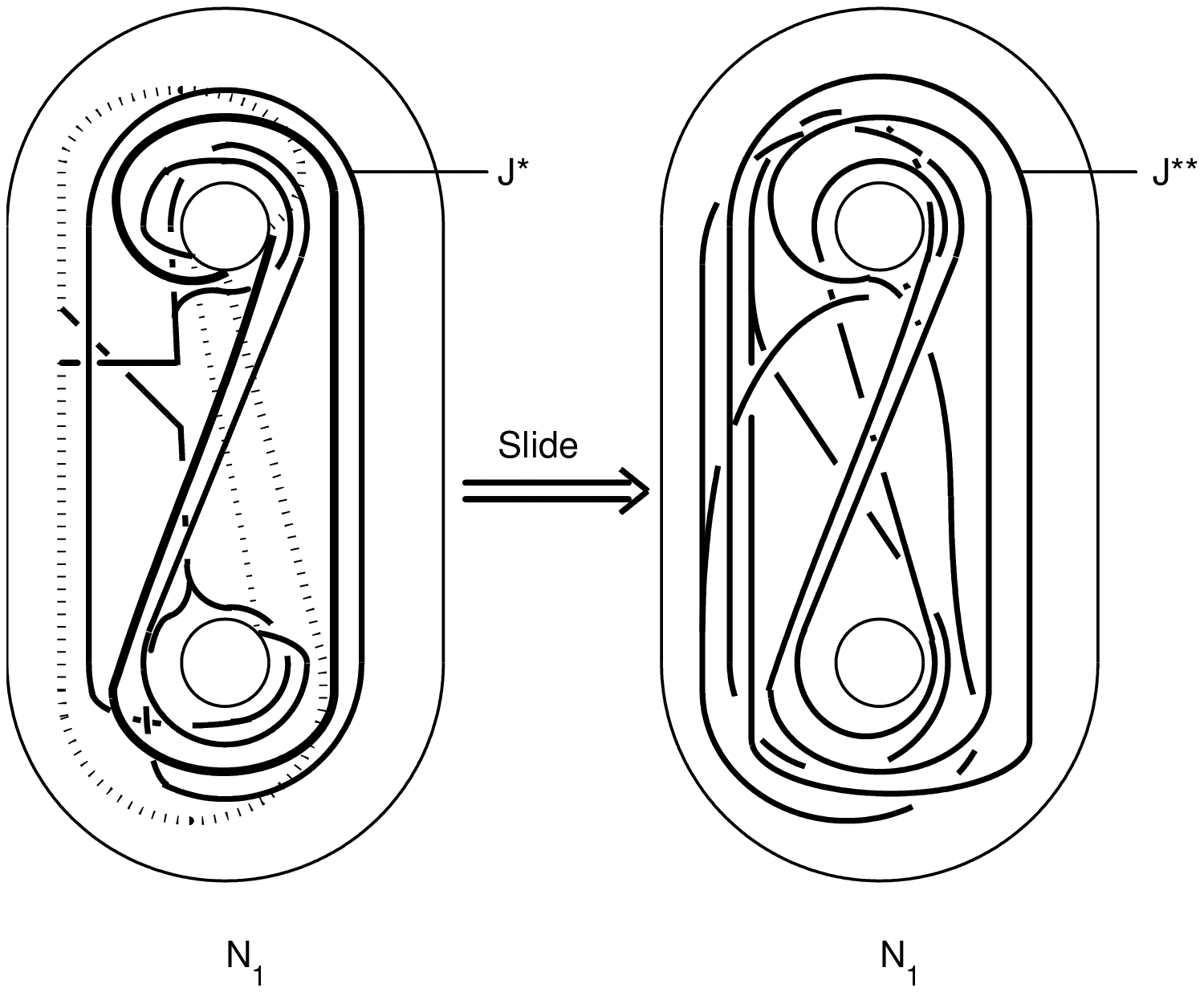}
\begin{center}
Figure 10
\end{center}
\end{center}

\begin{center}
\includegraphics[totalheight=6cm]{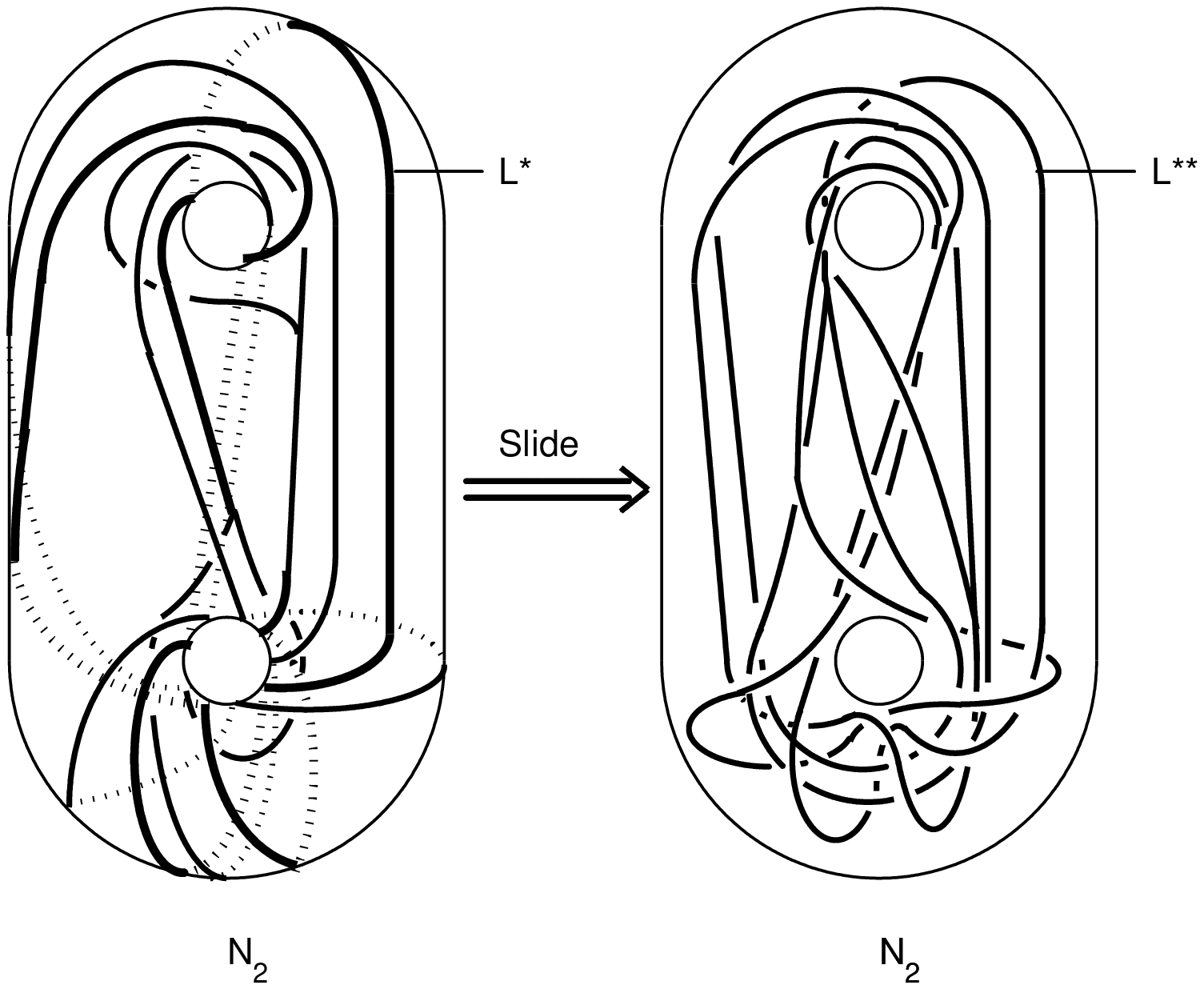}
\begin{center}
Figure 11
\end{center}
\end{center}

Now, we construct $f$. As Figures 10 and 11 show, we get $J^{**}$
($L^{**}$) from $J^{*}$ ($L^{*}$) by slide operations. Since
$J^{**}\sqcup L$ is isotopic to $L^{**}\sqcup J$ and as in the proof
of Proposition 4.3, we can construct $f\in Diff(M)$ which is
isotopic to identity, such that $f(J)=J^{**}$, $f(L^{**})=L$ and
$\Omega(f)$ is composed of type I Smale-Williams solenoids. Note
that here the alternating condition is used to show that the induced
matrix is irreducible, so Axiom 2 of
Smale-Williams solenoid follows.\\

\textbf{Theorem 4.5}:
\begin{itshape}If a Heegaard splitting  $M=N_1\cup N_2$ of the closed orientable
3-manifolds $M$ is a genus two alternating Heegaard splitting,
then there is a diffeomorphism $f$, such that $\Omega(f)$ consists
of two Smale-Williams solenoids.
\end{itshape}\\

\textbf{Proof}: It is the same as Example 4.4. Q.E.D.\\

\textbf{Example 4.6(The Truncated-Cube Space)}: Figure 12 is an
 alternating Heegaard diagram of a closed 3-manifold $M=N_1\cup N_2$ see from outside of $N_1$.

 %And Figure 13 is the same Heegaard splitting see from outside of $N_2$, it is from the
% Figure 14 by  Dehn twists $D_{e_3} D_{e_3} D_{e_2} D_{e_1}$  and then follows by a mirror
% symmetry.

 Its fundamental group,
$\pi_1(M)=<x_1,x_2;x_1x_2x_1x_2^{-1}x_1^{-1}x_2^{-1},x_1x_2x_1^{-1}x_2x_1x_2^{-1}
>$. Let $a=x_2$ and $b=x_2x_1$, we have
$\pi_1(M)=<a,b;a^{4}=b^{3}=(ab)^{2}>$, it is the extended triangle
group, $|\pi_1(M)|=48$, so $M$ is a genus two Seifert manifold with
base surface $S^2$ and three singular fibers, $M$ is called the
truncated-cube space, see [M] and [T].

\begin{center}
\includegraphics[totalheight=5cm]{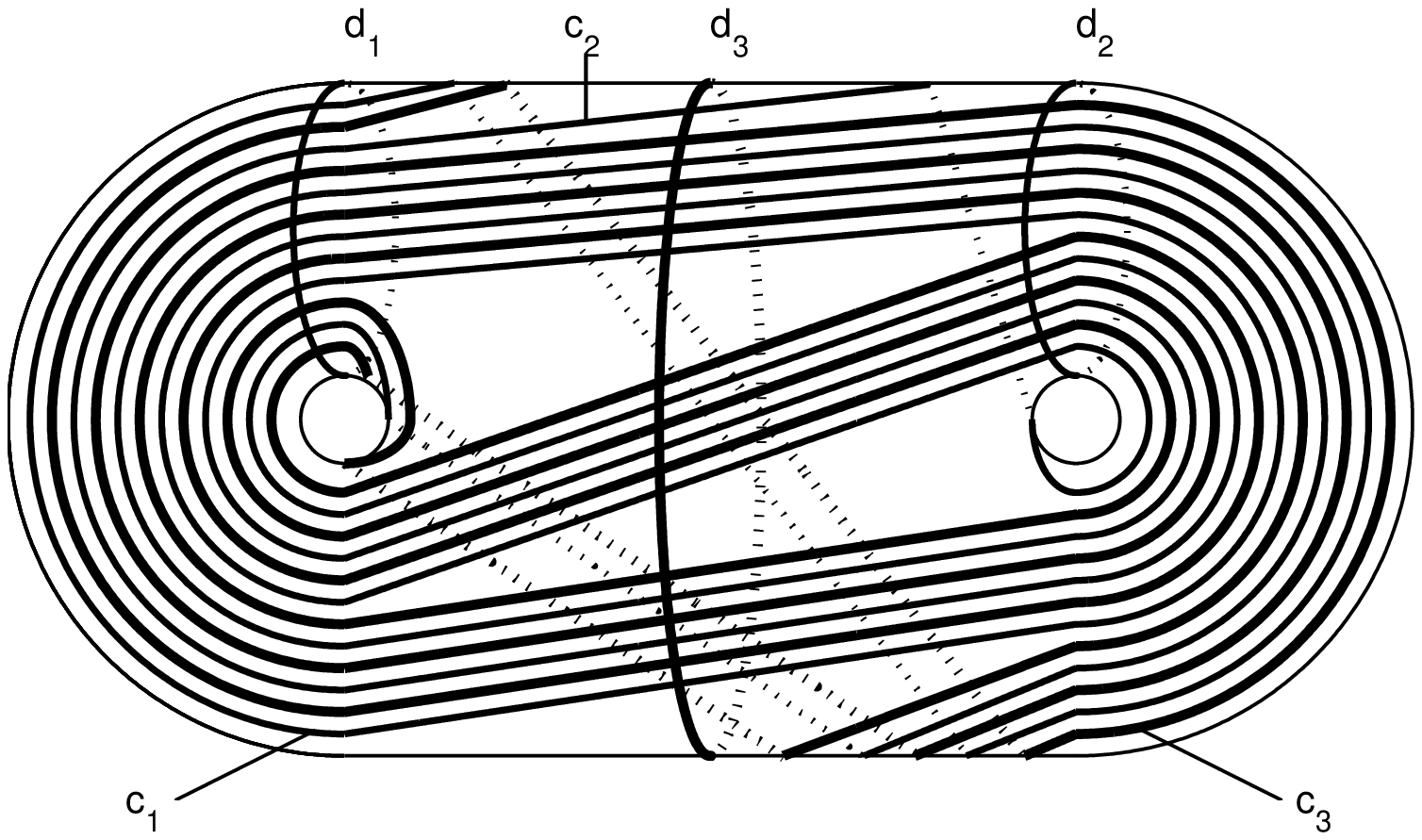}
\begin{center}
Figure 12
\end{center}
\end{center}

%\begin{center}
%\includegraphics[totalheight=6cm]{g2sws13.eps}
%\begin{center}
%Figure 13
%\end{center}
%\end{center}

%\begin{center}
%\includegraphics[totalheight=5cm]{g2sws14.eps}
%\begin{center}
%Figure 14
%\end{center}
%\end{center}

\textbf{Example 4.7($S^{3}$ with type II Smale-Williams solenoids)}:
Figure 13 is an alternating Heegaard splitting of $S^{3}$, so type
II Smale-Williams solenoids can be realized in $S^{3}$.

\begin{center}
\includegraphics[totalheight=6cm]{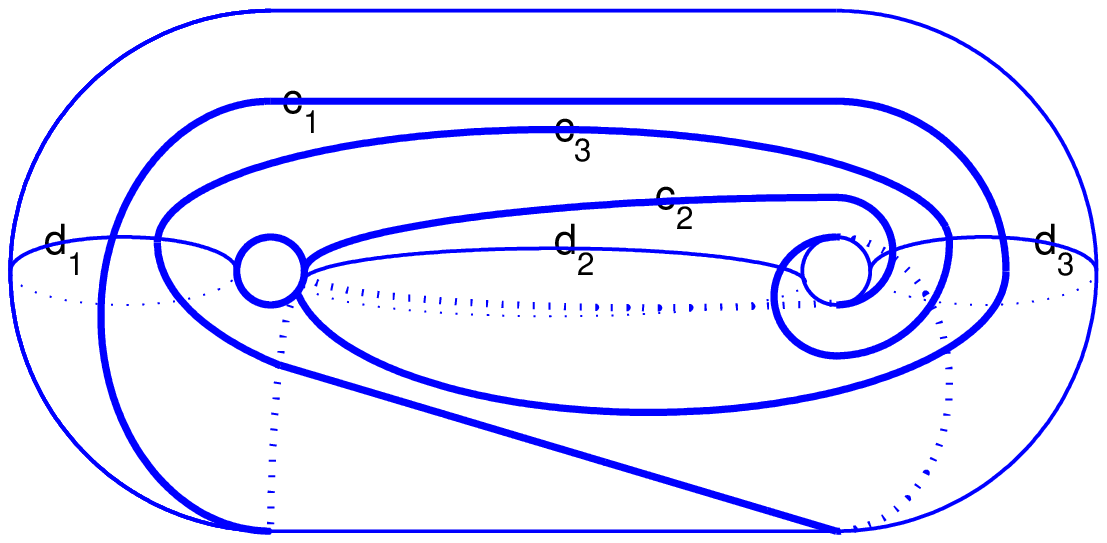}
\begin{center}
Figure 13
\end{center}
\end{center}

{\bf Acknowledgements}: The authors would like to thank Professor
Shicheng Wang for his many helpful discussions. Some of this work
was completed while the first named author was a postdoctoral
fellow at Peking University, and he thanks
the department for its hospitality.\\

{\bf References.} \vskip 4mm

[B1]  H. G. Bothe,  \emph{The ambient structure of expanding
attractors. I. Local triviality, tubular neighborhoods.} Math.
Nachr. 107 (1982), 327-348.

[B2]  H. G. Bothe, \emph{The ambient structure of expanding
attractors. II. Solenoids in 3-manifolds.} Math. Nachr. 112 (1983),
69-102.

[BH]  M. Bestvita and M. Handel, \emph{Train tracks and
automorphisms of
 free groups}, Ann. of Math. 135(1992), 1-51.

[G]  J. C. Gibbons, \emph{One-dimensional basic sets in the
three-sphere.}
 Trans. Amer. Math. Soc. 164 (1972), 163--178.

[H]  J. Hempel,\emph{ 3-Manifolds}, Ann. of Math. Studies, vol. 86,
1976.

[J]  W. Jaco,  \emph{Lectures on three-manifold Topology}, Published
by Amer. Math. Soc. 1980.

[JNW]  Boju Jiang, Yi Ni and Shicheng Wang, \emph{3-manifolds that
admit knotted solenoids as attractors},
 Trans. Amer. Math. Soc. 356, 4371-4382,(2004).

[M] M. Montesinos, \emph{Classical Tessellations and
three-manifolds}, Springer-Verlag, 1985.

[MY]  Jiming Ma, and Bin Yu, \emph{The realization of Smale solenoid
type attractors in 3-manifolds}, Topology and its
application,154(2007), no.17, 3021-3031.

[S]  S. Smale,\emph{ Differentiable dynamical systems}, Bull. Amer.
Math. Soc. 73 (1967), 747-817.

[T]  C. B. Thomas, \emph{Elliptic structures on 3-manifolds},
Cambridge University Press, 1986.

[W1]  R. F. Williams, \emph{One-dimensional non-wandering sets},
Topology Vol.6, pp.473-487, (1967).

[W2]  R. F. Williams, \emph{Expanding attractors}, Publ. Math. I. H.
E. S. 43,169-203, (1974).

%[W] YingQing Wu, \emph{Incompressible surfaces and Dehn surgery on
%$1$-bridge knots in handlebodies}.
% Math. Proc. Cambridge Philos. Soc. 120 (1996), no. 4, 687--696.

School of Mathematical sciences, Fudan University,  Shanghai
200433, P.R.China

 majiming@fudan.edu.cn

Department of Mathematics, Tongji University,  Shanghai 200092,
P.R.China

binyu1980@gmail.com

\end{document}